\documentclass[12pt]{article}
\usepackage{amssymb}
\begin{document}

\renewcommand{\thefootnote}{\arabic{footnote}}
\newcommand{\qqed}{\Box}
\newcommand{\blambda}{\mbox{\boldmath {$\mu$}}}
\newcommand{\bgamma}{\mbox{\boldmath {$\gamma$}}}
\newcommand{\bea}{\begin{eqnarray}}
\newcommand{\ena}{\end{eqnarray}}
\newcommand{\beas}{\begin{eqnarray*}}
\newcommand{\enas}{\end{eqnarray*}}
\newtheorem{theorem}{Theorem}[section]
\newtheorem{corollary}{Corollary}[section]
\newtheorem{conjecture}{Conjecture}[section]
\newtheorem{proposition}{Proposition}[section]
\newtheorem{lemma}{Lemma}[section]
\newtheorem{definition}{Definition}[section]
\newtheorem{example}{Example}[section]
\newtheorem{remark}{Remark}[section]

\def\blackslug{\hbox{\hskip 1pt \vrule width 4pt height 8pt depth 1.5pt
\hskip 1pt}}
\def\QED{\quad\blackslug\lower 8.5pt\null\par}
\def\proof{\par\penalty-1000\vskip .5 pt\noindent{\bf Proof\/: }}
\def\Var{{\rm Var}}

\begin{center}

{\Large \bf Stein's Method and the Zero Bias Transformation with Application to Simple Random Sampling}\end{center}
\vspace{3mm}

\begin{center}{Larry Goldstein and Gesine Reinert}\\
\end{center}
\footnotetext{Partially supported by NSF grant DMS-9505075}

\abstract Let $W$ be
a random variable with mean zero and variance $\sigma^2$.
The distribution of a variate $W^*$,
satisfying $EWf(W)=\sigma ^2 Ef'(W^*)$ for smooth functions $f$,
exists uniquely and defines the zero bias transformation
on the distribution of $W$.
The zero bias transformation shares many interesting properties with the
well known size bias transformation for non-negative variables,
but is applied to variables taking on both positive and negative values.
The transformation can also be defined on more general random objects.
The relation between the transformation and the expression
$wf'(w)-\sigma^2 f''(w)$ which appears in the Stein equation characterizing
the mean zero, variance $\sigma ^2$ normal $\sigma Z$ can be used to obtain
bounds on the difference $E\{h(W/\sigma)-h(Z)\}$ for smooth functions $h$ by constructing the pair $(W,W^*)$ jointly on the same space.  When $W$
is a sum of $n$ not necessarily independent variates, under certain conditions
which include a vanishing third moment, bounds on this difference
of the order $1/n$ for classes of smooth functions $h$ may be obtained.
The technique is illustrated by an application to simple random sampling.

\section{Introduction}
\label{intro}
Since 1972, Stein's method [\ref{st72}] has been extended and
refined by many authors and has become a valuable tool for deriving bounds for distributional approximations, in particular, for normal and Poisson approximations for sums of random variables. (In the normal case, see, for example, Ho and Chen [\ref{hoch78}], Stein [\ref{st86}], [\ref{st92}], Barbour [\ref{ADB90}], G\"otze [\ref{g91}], Bolthausen and  G\"otze [\ref{BG93}], Rinott [\ref{yos94}], and Goldstein and Rinott [\ref{gr96}]). Through the use of differential or difference equations which characterize the target distribution, Stein's method allows many different types of dependence structures to be
treated, and yields computable bounds on the approximation error.

The Stein equation for the normal is motivated by the fact that $W \sim {\cal N}(\mu,\sigma^2)$ if and only if
\beas
E\left\{(W-\mu)f'(W)-\sigma^2 f ''(W)\right\} =0 \quad \mbox{for all smooth $f$}.
\enas
Given a test function $h$, let $\Phi h=Eh(Z)$ where $Z \sim {\cal N}(0, 1)$.
If $W$ is close to $ {\cal N}(\mu, \sigma^2)$, $Eh((W-\mu)/\sigma)-\Phi h$ will be close to zero for
a large class of functions $h$, and $E\left\{(W-\mu)f'(W)-\sigma^2 f ''(W)\right\}$ will be
close to zero for a large class of functions $f$. It is natural then, given
$h$, to relate the functions $h$ and $f$ through the differential equation
\begin{equation}
\label{eq:ste}
(x-\mu)f'(x)-\sigma^2f''(x)=h((x-\mu)/\sigma)-\Phi h,
\end{equation}
and upon solving for $f$, compute $Eh((W-\mu)/\sigma)-\Phi h$ by $E\left\{(W-\mu)f'(W)-\sigma^2 f ''(W)\right\}$ for this $f$.
A bound on $Eh((W-\mu)/\sigma)-\Phi h$ can then be obtained by bounding the difference between  $E(W-\mu)f'(W)$ and $\sigma^2 Ef''(W)$.

Stein [\ref{st92}],  Baldi, Rinott, and Stein [\ref{brs}], and Goldstein and Rinott [\ref{gr96}], among others, were able to exploit a connection between
the Stein equation (\ref{eq:ste}), and the size biasing of nonnegative random variables.
If $W \ge 0$ has mean $0<EW=\mu < \infty,$
we say $W^s$ has the $W$-size biased distribution if for all
$f$ such that $EWf(W)$ exists,
\begin{equation}
\label{sbch}
EWf(W)=\mu Ef(W^s).
\end{equation}
The connection between the Stein equation and size biasing is described in Goldstein and  Rinott [\ref{gr96}]. In brief, one can obtain a bound on  $Eh((W-\mu)/\sigma)-\Phi h$ in terms of a pair $(W, W^s)$, coupled on a joint space, where $W^s$ has the $W$-size biased distribution. The terms in this bound will be small if $W$ and $W^s$ are close. The variates
$W$ and $W^s$ will be close,
for example, when $W=X_1+\cdots + X_n$ is the sum of i.i.d. random variables, as then $W^s$ can be constructed by replacing a single summand $X_i$ by an independent variate $X_i^s$ that has the $X_i$-size biased distribution. Similar constructions exist for
non-identically distributed and possibly dependent variates, and are studied in [\ref{gr96}].

As noted in [\ref{gr96}], the size biasing method  works well for combinatorial problems
such as counting the number of vertices in a random graph having prespecified degrees.  When the distributions approximated are counts, size biasing is natural; in particular, the counts $W$ are necessarily nonnegative. To size bias a $W$ which may take on both positive and negative values, it may be that for some $\rho$,  $W+\rho$ or $-W+\rho$ is a nonnegative random variable whose mean exists. Yet if $W$ has support on both the infinite positive and negative half lines then some truncation must be involved
in order to obtain a nonnegative random variable on which the size bias transformation can be performed. This is
especially unnatural if $W$ is symmetric, as one would expect that $W$
itself would be closer to normal than any version of itself involving translation and truncation.

The transformation and associated coupling which we study here has many similarities to the size biasing approach, yet it may be applied directly to mean zero random variables  and is particularly useful for symmetric random variables or those with vanishing third moment. The transformation is motivated by the size bias transformation and the equation that characterizes the mean zero normal:
\bea \label{zeromeanequ}
Z \sim {\cal N}(0,\sigma ^2) \quad  \mbox{if and only if} \quad  EWf(W)=\sigma^2 Ef'(W).
\ena
The similarity of the latter equation to equation (\ref{sbch})
suggests, given a mean
zero random variable $W$, considering a new distribution related to the distribution of $W$ according to the following definition.
\begin{definition}
\label{chzb}
Let $W$ be a mean zero random variable with finite, nonzero variance $\sigma^2$.
We say that $W^*$ has the $W$-zero biased distribution if for all differentiable $f$ for which $EWf(W)$ exists,
\begin{equation}
\label{cheq}
EWf(W)=\sigma ^2 Ef'(W^*).
\end{equation}
\end{definition}

The existence of the zero bias distribution for any such $W$
is easily established. For
a given $g \in C_c$, the collection of continuous functions with compact support, let $G=\int_0 ^w g$. The quantity
$$Tg = \sigma ^{-2} E\left\{WG(W)\right\}$$
exists since $EW^2 < \infty$, and defines a linear
operator $T:C_c \rightarrow {\bf R}$. To see moreover that $T$ is
positive, take $g \geq 0$. Then
$G$ is increasing, and therefore $W$ and $G(W)$ are positively
correlated. Hence $EWG(W) \geq EW EG(W) =0$, and $T$ is positive.
Now invoking the Riesz
representation theorem (see, eg. [\ref{foll}]),
we have
$Tg=\int g d\nu$ for some unique Radon measure $\nu$, which is
a probability measure by virtue of $T1=1$.  In fact,
the $W$-zero biased distribution is continuous for any nontrivial $W$; the density of $W^*$ is calculated explicitly in
Lemma \ref{zerokey}, part (\ref{2}).

Definition \ref{chzb} describes a transformation, which we term the zero bias transformation, on distribution functions with mean
zero and finite variance. However, for any $W$ with finite variance we can apply the transformation to the centered variate $W - EW$.

The zero bias transformation has many interesting properties, some of
which we collect below in Lemma \ref{zerokey}. In particular,
the mean zero normal is the unique fixed
point of the zero bias transformation. From this it is intuitive that
$W$ will be close to normal in distribution if and only if $W$ is close in distribution to $W^*$.

Use of the zero bias method, as with other like techniques, is through the use
of coupling and a Taylor expansion
of the Stein equation; in particular, we have
$$ E[Wf'(W) - \sigma^2 f''(W)] = \sigma^2 E[f''(W^*)-f''(W)],$$
and the right hand side may now immediately be expanded about $W$. In contrast, the use of other techniques such as size biasing requires an intermediate step which generates an additional error term (e.g., see equation (19) in  [\ref{gr96}]). For this reason, using the zero bias technique one is able to show why bounds of smaller order than $1/{\sqrt n}$ for smooth functions $h$ may be obtained when certain additional moment conditions apply.

For distributions with smooth densities, Edgeworth expansions reveal a similar phenomenon to what is studied here. For example, (see Feller [\ref{Feller2}]), if $F$ has a density and vanishing third moment, then an i.i.d. sum
of variates with distribution $F$ has a density which can be uniformly approximated by the normal to within a factor of $1/n$. However, these results depend on the smoothness of the parent distribution $F$. What we show here,
in the i.i.d. case say, is that for smooth test functions $h$,
bounds of order $1/n$ hold for any $F$ with vanishing third moment and finite fourth moment (see Corollary \ref{iid}).

Generally, bounds for non-smooth functions are more
informative than bounds for smooth functions (see for instance G\"otze [\ref{g91}], Bolthausen and G\"otze [\ref{BG93}], Rinott and Rotar [\ref{rinrot}] and Dembo and Rinott [\ref{yosefication}]); bounds
for non-smooth functions
can be used for the construction of confidence intervals, for instance.
Although the zero bias method can be used to obtain
bounds for non-smooth functions, we will
consider only smooth functions for the following reason.
At present, constructions for use
of the zero bias method are somewhat more difficult to achieve
than constructions for other methods; in particular,
compare the size biased construction in Lemma 2.1 of [\ref{gr96}] to
the construction in Theorem \ref{construction} here. Hence, for
non-smooth functions, other techniques may be easier to apply.
However, under added assumptions, the extra
effort in applying the zero bias method will be rewarded by
improved error bounds which may not
hold over the class of non-smooth functions. For example, consider
the i.i.d. sum of symmetric $+1,-1$ variates; the bound
on non-smooth functions of order $1/{\sqrt n}$
is unimprovable and may be obtained by a variety of methods.
Yet a bound of order $1/n$ holds for
smooth functions, and can be shown to be achieved by the
zero bias method. Hence, in order to reap the improved
error bound benefit of the zero bias method when such
can be achieved, we restrict attention to the class
of smooth functions.

Ideas related to the zero bias transformation have been studied by Ho and Chen [\ref{hoch78}], and Cacoullos et al. [\ref{cac2}]. Ho and Chen consider the zero bias distribution implicitly
(see equation 1.3 of [\ref{hoch78}]) in their version of one of Stein's proofs of the Berry Esseen theorem.
They treat a case with a $W$ the sum of dependent variates, and obtain rates of $1/ \sqrt{n}$ for the $L_p$ norm of the difference between the distribution function of $W$ and
the normal.

The approach of Cacoullos et al. [\ref{cac2}] is also related to what is studied here.
In the zero bias transformation, the distribution of $W$ is changed to that
of $W^*$ on the right hand side of identity (\ref{zeromeanequ}), keeping the form
of this identity, yielding (\ref{cheq}). In [\ref{cac2}], the distribution of $W$ is
preserved on the right hand side of (\ref{zeromeanequ}), and the form of the identity
changed to $E[Wf(W)]=\sigma ^2 E[u(W)f'(W)]$, with the function $u$ determined by the
distribution of $W$. Note that both approaches reduce to identity (\ref{zeromeanequ}) when $W$
is normal; in the first case $W^* \stackrel{d}{=} W$, and in the second, $u(w)=1$.

 The paper is organized as follows. In Section 2, we present some of the properties of the zero bias transformation and give two coupling constructions that generate $W$ and $W^*$ on a joint space. The first construction,
Lemma \ref{zerokey}, part \ref{5}, is
for the sum of independent variates, and its generalization,
Theorem \ref{construction},
for possibly dependent variates. In Section 3, we show how the zero bias transformation may be used to obtain bounds on the accuracy of
 the normal approximation in general. In Section 4, we apply the preceding results to obtain bounds of the order $1/n$ for smooth functions
$h$ when $W$ is
a sum obtained from simple random sampling without replacement (a case of global dependence), under a vanishing third moment assumption.

\section{The Zero Bias Transformation}

The following lemma summarizes some of the important features of the zero bias
transformation; property (\ref{lem:mom}) for $n=1$ will be
of special importance, as it gives that $EW^*=0$ whenever $EW^3=0$.

\begin{lemma} \label{zerokey}
Let $W$ be a mean zero variable with finite, nonzero variance, and
let $W^*$ have
the $W$-zero biased distribution in accordance with Definition \ref{chzb}.
Then,
\begin{enumerate}
\item  \label{1} The mean zero normal is the unique
fixed point of the zero bias transformation.

\item \label{2} The zero bias distribution is unimodal about zero and
continuous with density
function $p(w)=\sigma ^{-2}E[W,W>w]$. It follows that
the support of $W^*$ is the closed convex hull of the support of $W$
and that $W^*$ is bounded whenever $W$ is bounded.

\item The zero bias transformation preserves symmetry.

\item \label{lem:mom} $\sigma ^2E(W^*)^{n}=EW^{n+2}/(n+1)$ for $n \ge 1$.

\item \label{5} Let $X_1,\ldots,X_n$ be independent mean zero random
variables with $EX_i^2=\sigma_i^2$. Set $W=X_1+\cdots+X_n$, and $EW^2=\sigma^2$.  Let $I$ be a random index
independent of the $X's$ such that
\[ P(I=i)=\sigma_i^2/\sigma^2. \]
Let
$$W_i =W-X_i= \sum_{j \neq i} X_j.$$
Then $W_I+X_I^*$ has the $W$-zero biased distribution.
(This is analogous to size biasing a sum of non-negative
independent variates by replacing a variate chosen
proportional to its expectation by one chosen independently from
its size biased distribution; see Lemma 2.1 in [\ref{gr96}]).

\item \label{6} Let $X$ be mean zero with variance $\sigma_X^2$ and distribution $dF$. Let $(\hat{X}',\hat{X}'')$ have distribution
$$
d\hat{F}({\hat x}',{\hat x}'')=\frac{({\hat x}'-{\hat x}'')^2}{2 \sigma_X^2}dF(\hat{x}')dF(\hat{x}'').
$$
Then, with $U$ an independent uniform variate on $[0,1]$,
$U \hat{X'} + (1-U) \hat{X''}$
has the $X$-zero biased distribution.
\end{enumerate}
\end{lemma}

{\bf Proof of claims:}
\begin{enumerate}

\item This is immediate from Definition \ref{chzb} and
the characterization (\ref{zeromeanequ}).

\item The function $p(w)$ is increasing for $w<0$, and decreasing for $w>0$. Since $EW=0$,
$p(w)$ has limit 0 at both plus and minus infinity, and $p(w)$ must therefore
be nonnegative and unimodal about zero.
That $p$ integrates to 1 and is the density of a variate $W^*$
which satisfies (\ref{cheq}) follows by uniqueness (see the remarks
following Definition \ref{chzb}), and by applying Fubini's theorem separately
to $E[f'(W^*);W^* \ge 0]$ and $E[f'(W^*);W^* < 0]$, using
\[
E[W;W>w]=-E[W;W\le w],
\]
which follows from $EW=0$.

\item If $w$ is a continuity point of the distribution function
of a symmetric $W$, then $E[W;W>w]=E[-W;-W>w]=-E[W;W<-w]=E[W;W>-w]$
using $EW=0$.
Thus, there is a version of the $dw$ density of $W^*$ which is the same
at $w$ and $-w$ for almost all $w$ [$dw$]; hence $W^*$ is symmetric.

\item Substitute $w^{n+1}/(n+1)$ for $f(w)$ in the characterizing equation
(\ref{cheq}).

\item Using independence and equation (\ref{cheq}) with $X_i$ replacing $W$,
\beas
\sigma^2 Ef'(W^*)&=&EWf(W)\\
&=&\sum_{i=1}^n EX_i f(W)\\
&=&\sum_{i=1}^n EX_i^2 Ef'(W_i+X_i^*)\\
&=&\sigma^2 \sum_{i=1} ^n \frac{\sigma_i ^2 }{\sigma^2} Ef'(W_i+X_i^*)\\
&=& \sigma^2 Ef'(W_I+X_I^*).\\
\enas
Hence, for all smooth $f$, $Ef'(W^*)=Ef'(W_I+X_I^*)$, and the
result follows.

\item Let $X',X''$ denote independent copies of the variate $X$. Then,
\beas
\sigma_X^2 Ef'(U{\hat X}'+(1-U){\hat X}'')
&=& \sigma_X^2 E\left(\frac{f({\hat X}')-
             f({\hat X}'')}{{\hat X}'-{\hat X}''}\right)\\
                       &=& \frac{1}{2} E(X'-X'')(f(X')-f(X''))\\
                       &=& EX'f(X')-EX''f(X')\\
                       &=& EXf(X)\\
                       &=&\sigma_X^2 Ef'(X^*).\\
\enas
Hence, for all smooth $f$,  $Ef'(U{\hat X}'+(1-U){\hat X}'')=Ef'(X^*)$.

\end{enumerate} $\hfill \qqed$

By (\ref{1}) of Lemma \ref{zerokey}, the mean zero normal is a fixed point of the zero bias transformation. One can also gain some insight
into the nature of the transformation by observing its action
on the distribution of the variate $X$
taking the values $-1$ and $+1$ with equal probability.
Calculating the density function of the $X$-zero
biased variate $X^*$ according to (\ref{2}) of Lemma \ref{zerokey},
we find that $X^*$ is uniformly distributed on the interval $[-1,1]$.
A similar calculation for the discrete mean zero variable $X$ taking
values $x_1 < x_2 < \cdots <x_n$ yields that the $X$-zero biased
distribution is a mixture of uniforms over the intervals $[x_i,x_{i+1}]$.
These examples may help in understanding how a uniform variate $U$
enters in (\ref{6}) of Lemma \ref{zerokey}.

For a construction of $W$ and $W^*$ which may be applied in the presence of dependence, in the remainder of this section, we will consider the following
framework. Let $X_1,\ldots,X_n$ be mean zero random variables, and with $W=X_1+\cdots+X_n$, suppose $EW^2=\sigma^2$ exists.
For each $i=1,\ldots,n$, assume that there exists a distribution $dF_{n,i}(x_1,\ldots,x_{i-1},x_i',x_i'',x_{i+1},\ldots,x_n)$ on $n+1$
variates $X_1,\ldots,X_{i-1},X_i',X_i'',X_{i+1},\ldots,X_n$
such that
\bea
\label{dieq1}
(X_1,\ldots,X_{i-1},X_i',X_i'',X_{i+1},\ldots,X_n) &\stackrel{d}{=}&
(X_1,\ldots,X_{i-1},X_i'',X_i',X_{i+1},\ldots,X_n), \nonumber \\
\ena
and

\bea
\label{dieq2}
(X_1,\ldots,X_{i-1},X_i,X_{i+1},\ldots,X_n) &\stackrel{d}{=}&
(X_1,\ldots,X_{i-1},X_i',X_{i+1},\ldots,X_n).
\ena
(The choice $X_i'=X_i$ is natural, and then (\ref{dieq2}) is satisfied.)

Further, we will suppose that there is a $\rho$ such that for all $f$ for which $EWf(W)$ exists,
\bea
\label{rho} \quad \sum_{i=1}^n EX_i' f(W_i+X_i'')&=&\rho EWf(W),
\ena
where $W_i = W-X_i$. We set
\bea \label{vi}
v_i^2=E(X_i'-X_i'')^2.
\ena
Under these conditions, we have the following proposition.

\begin{proposition} \label{rhoconstruction}
\beas
\rho = 1 - \frac{1}{2 \sigma^2} \sum_{i=1}^n v_i^2.
\enas
\end{proposition}

Before proving this proposition, note that if a collection of variates already satisfies
(\ref{dieq1}) and (\ref{dieq2}), and that if for each $i$,
\begin{equation}
\label{exhrho}
E\{X_i'|W_i+X_i''\}={\rho \over n} (W_i+X_i''),
\end{equation}
then
$$
EX_i'f(W_i+X_i'')=\frac{\rho}{n}EWf(W),
$$
and so condition (\ref{rho}) will be satisfied.

\noindent {\bf Proof of Proposition \ref{rhoconstruction}:}
Substituting $f(x) =x$ in (\ref{rho}) yields, by (\ref{dieq2}), that
\begin{eqnarray*}
\rho \sigma^2 &=& \sum_{i=1}^n EX_i'(W_i+X_i'')\\
&=& \sum_{i=1}^n EX_i(W - X_i) + \sum_{i=1}^n EX_i'X_i''\\
&=& \sigma^2 - \sum_{i=1}^n  EX_i^2 - {1 \over 2} \sum_{i=1}^n \{ E(X_i' - X_i'')^2 -  E (X_i')^2 - E(X_i'')^2\}\\
&=&  \sigma^2  - {1 \over 2} \sum_{i=1}^n v_i^2,
\end{eqnarray*}
so that Proposition (\ref{rhoconstruction}) follows.
 $\qqed$

The following theorem, generalizing (\ref{5}) of
Lemma \ref{zerokey}, gives a coupling construction for
$W$ and $W^*$ which may be applied in the presence of dependence
under the framework of Proposition \ref{rhoconstruction}.

\begin{theorem} \label{construction}
Let $I$ be a random index
independent of the $X's$ such that
\[ P(I=i)=v_i^2/\sum_{j=1}^n v_j^2.\]
Further, for $i$ such that $v_i>0$, let $\hat{X}_1,\ldots,\hat{X}_{i-1},\hat{X}_i',\hat{X}_i'',\hat{X}_{i+1},\ldots,\hat{X}_n$
be chosen according to the distribution
\bea
\label{hatF}
\lefteqn{d\hat{F}_{n,i}(\hat{x}_1,\ldots,\hat{x}_{i-1},\hat{x}_i',\hat{x}_i'',\hat{x}_{i+1},\ldots,\hat{x}_n)} \nonumber \\
&=& \frac{(\hat{x}_i'-\hat{x}_i'')^2}{v_i^2}dF_{n,i}(\hat{x}_1,\ldots,\hat{x}_{i-1},\hat{x}_i',\hat{x}_i'',\hat{x}_{i+1},\ldots,\hat{x}_n).
\ena
Put
\begin{equation}
\label{Put}
 \hat{W}_i = \sum_{j \neq i} \hat{X}_j.
\end{equation}
Then, with $U$ a uniform $U[0,1]$ variate which is independent
of the $X$'s and the index $I$,
\[ U \hat{X}_I' + (1-U) \hat{X}_I''+ \hat{W}_I  \]
has the $W$-zero biased distribution.

In particular, when $X_1,\ldots,X_n$ are exchangeable, if one constructs
exchangeable variables with distribution $dF_{n,1}$ which satisfy $v_1^2>0$,
(\ref{dieq2}), and (\ref{exhrho}) for $i=1$,
then
$$
U \hat{X}_1' + (1-U) \hat{X}_1''+ \hat{W}_1
$$
has the $W$-zero biased distribution.
\end{theorem}


\noindent {\bf Proof of Theorem \ref{construction}:}
With Proposition \ref{rhoconstruction} we have
\beas
E\int_0^1 f'(u\hat{X}_I' + (1-u) \hat{X}_I''+ \hat{W}_I) du
&=& E\left( \frac{f(\hat{W}_I  + \hat{X}_I') - f(\hat{W}_I + \hat{X}_I'')}{\hat{X}_I' - \hat{X}_I''}\right) \\
&=& \sum_{i=1}^n \frac{v_i^2}{\sum v_j^2}   E\left( \frac{f(\hat{W}_i + \hat{X}_i') - f(\hat{W}_i + \hat{X}_i'')} {\hat{X}_i' - \hat{X}_i''}\right) \\
&=& \frac{1}{\sum v_j^2} \sum_{i=1}^n   E(X_i'-X_i'')(f(W_i + X_i') -
f(W_i + X_i'')) \\
&=&\frac{2}{\sum v_j^2} \sum_{i=1}^n \left( E X_i' f(W_i + X_i') - E X_i' f(W_i + X_i'') \right)\\
&=& \frac{2}{\sum v_j^2} \left\{ EWf(W)-\rho EWf(W) \right\} \\
&=& \frac{2(1-\rho)}{\sum v_j^2} EWf(W)\\
&=& \frac{1}{\sigma^2} EWf(W)\\
&=& Ef'(W^*), \enas
using Proposition (\ref{rhoconstruction}) for the next to last step.

To show the claim
in the case where the variates are exchangeable, set $dF_{n,i}=dF_{n,1}$
for $i=2,\ldots,n$ and observe that the $dF_{n,i}$ so
defined now satisfy the conditions of the theorem, and
the distributions of the resulting $U{\hat X}'_i+(1-U){\hat X}''_i+{\hat W}_i$
does not depend on $i$.

$\hfill \qqed$

Note that if the variates $X_1,\ldots,X_n$ are independent, one can generate the collection $X_1,\ldots,X_{i-1},X_i',X_i'',X_{i+1},\ldots,X_n$ by letting $X_i',X_i''$ be independent
replicates of $X_i$. In this case, conditions (\ref{dieq1}), (\ref{dieq2}), and (\ref{rho}) above are
satisfied, the last with $\rho=0$, and
the construction reduces to that given in (\ref{5}) of Lemma \ref{zerokey},
in view of (\ref{6}) of that same lemma.

\section{Bounds in the Central Limit Theorem}
The construction of Lemma \ref{zerokey}, part \ref{5},
together with the following bounds of Barbour [\ref{ADB90}] and G\"{o}tze  [\ref{g91}] on the solution $f$ of the differential equation (\ref{eq:ste}) for a test function $h$ with $k$ bounded derivatives,
\begin{equation}
\label{solnbds}
||f^{(j)}|| \le ({j}\sigma^{j})^{-1} ||h^{(j)}||\quad \mbox{$j=1,\ldots,k$,}
\end{equation}
yield the following remarkably simple proof of the Central Limit Theorem, with
bounds on the approximation error, for independent possibly non-identically
distributed mean zero variables $X_1,\ldots,X_n$
with variance 1 and common absolute first and third moments.

By Lemma \ref{zerokey}, part (\ref{5}), using independence, we can achieve $W^*$ having the $W$-zero biased distribution by selecting a random index $I$ uniformly and replacing $X_I$ by an independent variable $X_I^*$ having the $X_I$-zero biased distribution. Now, since
$EWf(W)=\sigma^2 Ef'(W^*)$, using the bound (\ref{solnbds}),
\begin{eqnarray}
|E\left\{h(W/\sigma)-\Phi h\right\}|&=&\left|E\left\{Wf'(W)-\sigma^2 f''(W)\right\}\right| \nonumber \\
              &=&\sigma^2 |E\left\{f''(W^*)-f''(W)\right\}|  \nonumber\\
              &\le&\sigma^2 ||f^{(3)}|| E|W^*-W|\nonumber \\
              &\le&\frac{1}{3\sigma}||h^{(3)}||E|X_I^*-X_I| \label{stop}.
\end{eqnarray}
Now, using the bound $E|X_I^*-X_I|\le E|X_I^*|+E|X_I|$ and the function $x^2\mbox{sgn}(x)$ and its derivative $2|x|$ in equation (\ref{cheq}) , we derive
$E|X_i^*|=\frac{1}{2}E|X_i|^3$, and therefore $E|X_I^*|=\frac{1}{2} E|X_1|^3$. Next, $E|X_I|=E|X_1|$, and by
H\"older's inequality and $EX_i^2=1$, we have $E|X_i| \le 1 \le  E|X_i|^3$. Hence, since $EW^2=n=\sigma^2$,
\beas
|E\left\{h(W/\sigma)-\Phi h\right\} | \le \frac{||h^{(3)}|| E|X_1|^3}{2 \sqrt{n}}.
\enas
Thus  we can obtain a bound of order $n^{-1/2}$ for
smooth test functions with an
explicit constant using only the first term in the Taylor expansion of $f''(W^*)-f''(W)$. For arbitrary independent mean zero variates, continuing
from (\ref{stop}), for small additional effort the above inequality
generalizes to
$$
|E\left\{h(W/\sigma)-\Phi h\right\} | \le \frac{||h^{(3)}|| E|X_I|^3}{2 \sigma}.
$$

The following theorem shows how
the distance between an arbitrary mean zero, finite variance random variable $W$ and a mean zero normal with the same variance can be bounded by
the distance between $W$ and a variate $W^*$ with the $W$-zero biased
distribution defined on a joint space. It is instructive to compare the following
theorem with Theorem 1.1 of [\ref{gr96}], the corresponding result
when using the
size biased transformation.
\begin{theorem} \label{zerogen}  Let $W$ be a mean zero random variable
with variance $\sigma^2$, and suppose $(W, W^*)$ is given on a joint probability space so
that $W^*$ has the $W$-zero biased distribution. Then for all $h$ with four bounded derivatives,
$$
|Eh(W/\sigma)-\Phi h| \le \frac{1}{3  \sigma}\|h^{(3)}\| \sqrt{E\{E(W^* - W|W)^2\}} +
\frac{1}{8  \sigma^2} \|h^{(4)}\| E(W^*-W)^2.
$$
\end{theorem}

\noindent {\bf Proof.}
For the given $h$, let $f$ be the solution to (\ref{eq:ste}). Then, using the bounds in (\ref{solnbds}),
it suffices to prove
$$
|E[Wf'(W) - \sigma^2f''(W)]| \leq \sigma^2 \|f^{(3)}\| \sqrt{E\{E(W^* - W|W)^2\}} +
\frac{\sigma^2}{2} \|f^{(4)}\| E(W^*-W)^2.
$$
By Taylor expansion, we have
\beas
|E[Wf'(W)- \sigma^2f''(W)]| & = &|  \sigma^2 E[f''(W^*) - f''(W)]|\\
              & \leq & \sigma^2 | E  f^{(3)}(W)(W^*-W) | +
\frac{\sigma^2}{2}
\|f^{(4)}\| E (W^*-W)^2.
\enas
For the first term, condition on $W$ and then apply the Cauchy-Schwarz inequality;
\beas
\big|E\big[f^{(3)}(W)(W^*-W)\big]\big| &= &\big|E\big[f^{(3)}(W)E(W^*-W|W)\big]\big| \\
&\leq& \|f^{(3)}\| \sqrt{E\{E(W^* - W|W)^2\}}. \qqed
\enas

For illustration only, we apply Theorem \ref{zerogen} to the sum of independent identically distributed variates  to show how the the zero bias transformation leads to an error bound for smooth functions of order $1/n$, under additional
moment assumptions which include a vanishing third moment.
\begin{corollary} \label{iid}
Let $X,X_1,X_2,\ldots,X_n$ be independent and identically distributed mean zero, variance one random variables with vanishing third moment and $EX^4$ finite. Set $W=\sum_{i=1}^nX_i$. Then  for any
function $h$ with four bounded derivatives,
\[
|E\left\{h(W/{\sqrt n})-\Phi h\right\}|\le  n^{-1}\left\{\frac{1}{3}||h^{(3)}||+ \frac{1}{6}|| h^{(4)}|| EX^4 \right\}.
\]
\end{corollary}

\noindent {\bf Proof:} Construct $W^*$ as in Lemma \ref{zerokey}, part (\ref{5}). Then
$$
E(W^* - W|W) = E(X_I^*- X_I|W)= E(X_I^*)-E (X_I|W),
$$
since $X_I^*$ and $W$ are independent.
Using the moment relation $EX^*=(1/2)EX^3$ given in Lemma \ref{zerokey}, part (4),  $EX_i^3=0$ implies that $EX_i^*=0$, and so $EX_I^*=0$.
 Using that the $X$'s are i.i.d., and therefore exchangeable, $E(X_I|W)=W/n$. Hence we obtain
$E(X_I^*-X_I|W)=-W/n,$ and
$$\sqrt{E\{E(X_I^*-X_I|W)^2\}} = {1 \over {\sqrt n}}.$$
For the second term in Theorem \ref{zerogen},
$$ E(W^*-W)^2 = E(X_I^*-X_I)^2.$$
The moment relation property (\ref{lem:mom}) in Lemma
\ref{zerokey} and the assumption that $EX^4$ exists renders
$E(X_I^*-X_I)^2$ finite and equal to $EX^4/3 + EX^2\le (4/3)EX^4$,
by $EX^2=1$ and H{\"o}lder's inequality. Now using $\sigma^2=n$
and applying Theorem \ref{zerogen} yields the assertion. $\hfill
\qqed$

It is interesting to note that the constant $\rho$ of equation (\ref{rho})
does not appear in the bounds of Theorem \ref{zerogen}. One explanation
of this phenomenon is as follows. The $\rho$ of the coupling of
Theorem \ref{construction} is related to the $\lambda \in (0,1)$ of a
coupling of Stein [\ref{st86}], where a mean zero exchangeable pair
$(W, W')$, with distribution $dF(w,w')$, satisfies
$E\{W'|W\}=(1-\lambda)W$. One can show that if $(\hat{W},\hat{W}')$
has distribution
$$
d\hat{F}(\hat{w},\hat{w}')=\frac{(\hat{w}-\hat{w}')^2}{E(W-W')^2} dF(\hat{w},\hat{w}'),
$$
then with $U$ a uniform variate on [0,1], independent of
all other variables,
$U\hat{W}+(1-U)\hat{W}'$ has the $W$-zero bias distribution.
Taking simple cases, one can see that the value of $\lambda$
has no relation of the closeness of $W$ to the normal. For
instance, if $W$ is the sum of $n$ i.i.d. mean zero, variance
one variables, then $W$ is close to normal when $n$ is large.
However, for a given value of $n$, we may achieve any $\lambda$
of the form $j/n$ by taking $W'$ to be the sum of any $n-j$ variables
that make up the sum $W$, added to $j$ i.i.d. variables that are independent
of those that form $W$, but which have the same distribution.

We only study here the notion of zero biasing in one
dimension; it is possible to extend this concept to any finite
dimension. The definition of zero biasing in ${\bf R}^p$ is
motivated by the following multivariate characterization.
A vector ${\bf Z} \in R^p$ is multivariate normal
with mean zero and covariance matrix ${\bf \Sigma}=(\sigma_{ij})$
if and only if for
all smooth test functions $f:R^p \rightarrow R$,
$$
E\sum_{i=1} ^p Z_i f _i({\bf Z}) = E\sum_{i,j=1}^p \sigma_{ij}f_{ij}({\bf Z}),
$$
where $f_i, f_{ij},\ldots$ denote the partial
derivatives of $f$ with respect to the indicated coordinates.
Guided by this identity, given a mean zero
vector ${\bf X}= (X_1, \ldots, X_p)$ with covariance matrix $\Sigma$, we say
the collection of vectors ${\bf X}^*=({\bf X}_{ij}^*)$ has the ${\bf X}$-zero bias distribution if
\begin{equation} \label{mvzb}
E \sum_{i=1}^p X_i f_i({\bf X})=E\sum_{i,j=1} ^p \sigma _{ij} f_{ij}
({\bf X}_{ij}^*),
\end{equation}
for all smooth $f$. As in the univariate case, the mean zero normal is a fixed point
of the zero bias transformation; that is, if ${\bf X}$ is a mean
zero normal vector, one may satisfy (\ref{mvzb}) by setting
${\bf X}_{ij}^*={\bf X}$ for all $i,j$.

Using the definition of zero biasing in finite dimension, one
can define the zero bias concept for random variables
over an arbitrary index set ${\cal H}$
as follows. Given a collection $\{\xi (\phi), \phi \in {\cal H}\}$ of
real valued mean zero random variables with nontrivial finite second moment,
we say the collection $\{\xi_{\phi \psi}^*, \phi, \psi \in {\cal H}\}$ has the
$\xi$-zero biased distribution if for all $p=1,2,\ldots$ and $(\phi_1,\phi_2,\ldots,\phi_p)\in {\cal H}^p$,
the collection of $p$-vectors $({\bf X}_{ij}^*)$ has the
${\bf X}$-zero bias distribution, where, for $1 \leq i, j \leq p$,
$$
({\bf X}_{ij}^*)=(\xi^*_{\phi_i \phi_j}(\phi_1),\ldots,\xi^*_{\phi_i \phi_j}(\phi_p)),
$$
and
$$
{\bf X}=(\xi(\phi_1),\ldots,\xi(\phi_p)).
$$
Again when $\xi$ is normal, we may set $\xi_{\phi \psi}^*=\xi$ for all $\phi, \psi$. This definition reduces to the one given above for random vectors when ${\cal H}
=\{1,2,\ldots,n\}$, and can be applied to, say, random processes
by setting ${\cal H}={\bf R}$, or random measures by letting
${\cal H}$ be a specified class of functions.

\section{Application: Simple random sampling}

We now apply Theorem \ref{zerogen} to obtain a bound on the error
incurred when using the normal to approximate the distribution of
a sum obtained by simple random sampling.
In order to obtain a bound of order $1/n$ for smooth functions, we impose an additional moment condition as in Corollary \ref{iid}.

Let ${\cal A}=\{a_1, \ldots., a_N\}$ be a set of
real numbers such that
\bea \label{moment}
\sum\limits_{a \in {\cal A}}
a = \sum\limits_{a \in {\cal A}} a^3 = 0;
\ena
the following is a useful consequence of (\ref{moment}),
\begin{equation}
\label{sym}
\mbox{for any }E \subset \{1, \ldots, N\} \quad
\mbox{and } k \in \{1,3\}, \quad \sum_{a \in E}a^k=-\sum_{a \not \in E}a^k.
\end{equation}

We assume until the statement of Theorem \ref{srs} that the
elements of ${\cal A}$
are distinct; this condition will be dropped in the theorem.
Let $0<n<N$, and set $N_n = N(N-1)\cdots(N-n+1)$, the $n^{th}$ falling factorial of $N$.
Consider the random vector ${\bf X}=(X_1,\ldots, X_n)$ obtained by a
simple random sample of size $n$ from ${\cal A}$, that is, ${\bf X}$ is a realization
of one of the equally likely $N_n$ vectors of distinct elements of ${\cal A}$. Put
\begin{equation}
\label{wsum}
W= X_1 + \cdots + X_n.
\end{equation}
Then, simply we have $EX_i=EX_i^3=EW=EW^3=0$, and
\beas
EX_i^2 = \frac{1}{N} \sum_{a \in {\cal A}} a^2= \sigma_X^2, \quad
EW^2 = \frac{n(N-n)}{N(N-1)} \sum_{a \in {\cal A}} a^2 =   \sigma^2, \quad \mbox{say.}
\enas
As we will consider the normalized variate $W/\sigma$, without loss of generality
we may assume
\bea \label{<2>}
\sum_{a \in {\cal A}} a^2 =1;
\ena
note that (\ref{<2>}) can always be enforced by rescaling ${\cal A}$,
leaving (\ref{moment}) unchanged.

The next proposition shows how to apply Theorem \ref{construction} to construct $W^*$ in the context of simple random sampling.

\begin{proposition} \label{zerosrs}
Let
\bea
\label{srsnp1}
dF_{n,1}(x_1',x_1'',x_2,\ldots,x_n)=N_{n+1}^{-1}{\bf 1}(\{x_1',x_1'',x_2,\ldots, x_n \} \subset {\cal A},
\mbox{ distinct}),
\ena
the simple random sampling distribution on $n+1$ variates from ${\cal A}$, and
$\hat{\bf X}=(\hat{X}'_1, \hat{X}''_1, \hat{X}_2, \ldots,  \hat{X}_n)$ be a random
vector with distribution
\bea \label{(**)}
d\hat{F}_{n,1}(\hat{\bf x})= \frac{(\hat{x}_1'-\hat{x}_1'')^2}{2N} (N-2)_{n-1}^{-1} {\bf 1} ( \{ \hat{x}'_1,
\hat{x}''_1, \hat{x}_2, \ldots, \hat{x}_n\} \subset {\cal A}, \mbox{ distinct}).
\ena
Then, with $U$ a uniform $[0,1]$ random variable independent of $\hat{\bf X}$, and
$\hat{W}_1$ given by (\ref{Put}),
\begin{equation}
\label{wstsum}
W^*= U\hat{X}_1' + (1-U) \hat{X}_1'' + \hat{W}_1
\end{equation}
has the $W$-zero biased distribution.
\end{proposition}

\noindent {\bf Proof.} We apply Theorem \ref{construction} for exchangeable variates. With $X_1,\ldots,X_n$ a simple random sample of size $n$,
the distributional relation (\ref{dieq2}) is immediate.
Next, using  the scaling (\ref{<2>}), we see that
$v_1^2$ given in (\ref{vi}) equals $2/(N-1)$, which is positive, and that furthermore, the distribution
(\ref{(**)}) is constructed from the distribution
(\ref{srsnp1}) according to the prescription  (\ref{hatF}).
Lastly, using (\ref{sym}) with $k=1$, we have
\beas
E\{X_1'|X_1'',X_2,\ldots,,X_n \}=-\left({{W_1 + X_1''} \over {N-n}}\right),
\enas
and hence condition (\ref{exhrho}) is satisfied with
$\rho=-n/(N-n)$. $\hfill \qqed$

We now begin to apply
Theorem \ref{zerogen} by constructing $W$ and $W^*$ on a joint space.
We achieve this goal by constructing the simple random sample ${\bf X}=(X_1,\ldots,X_n)$ together with the variates
$\hat{\bf X}=(\hat{X}'_1, \hat{X}''_1, \hat{X}_2, \ldots, \hat{X}_n)$ with distribution as in (\ref{(**)}) of Proposition \ref{zerosrs}; $W$ and $W^*$
are then formed from these
variates according to (\ref{wsum}) and (\ref{wstsum}) respectively.

{\bf Construction of $W$ and $W^*$.}
Start the construction with the simple random sample ${\bf X}=(X_1,\ldots,X_n)$.
To begin the construction of $\hat{\bf X}$ with distribution (\ref{(**)}),
set
$$
q(u,v) = \frac{(u-v)^2}{2 N}{\bf 1}(\{u,v \} \subset {\cal A}).
$$
Note that variates $U,V$ with distribution $q(u,v)$ will be unequal, and
therefore we have that the distribution (\ref{(**)}) factors as
\begin{eqnarray} \label{factor}
d\hat{F}_{n,1}(\hat{\bf x})= q(\hat{x}_1',\hat{x}_1'')(N-2)_{n-1}^{-1} {\bf 1} ( \{ \hat{x}_2, \ldots, \hat{x}_n\} \subset {\cal A}\setminus \{\hat{x}'_1,
\hat{x}''_1\} , \mbox{ distinct}).
\end{eqnarray}
Hence, given $(\hat{X}'_1, \hat{X}''_1)$, the vector $(\hat{X}_2,\ldots,\hat{X}_n)$
is a simple random sample of size $n-1$ from the $N-2$ elements
of ${\cal A}\setminus \{\hat{X}'_1,\hat{X}''_1\}$.

Now, independently of the chosen sample ${\bf X}$, pick
$(\hat{X}'_1, \hat{X}''_1)$ from the distribution $q(u,v)$.
The variates $(\hat{X}'_1, \hat{X}''_1)$ are then placed as the first two components in the
vector $\hat{\bf X}$. How the remaining $n-1$ variates in $\hat{\bf X}$ are chosen depends
on the amount of intersection between the sets $\{X_2,\ldots,X_n\}$ and
$\{\hat{X}'_1, \hat{X}''_1\}$. If these two sets do not intersect, fill in the remaining $n-1$
components of $\hat{\bf X}$ with $(X_2,\ldots,X_n)$. If the sets have an intersection,
remove from the vector $(X_2,\ldots,X_n)$ the two variates (or single variate)
that intersect and replace
them (or it) with values obtained by a simple random sample of size two (one) from
${\cal A}\setminus \{\hat{X}'_1, \hat{X}''_1, X_2, \ldots, X_n\}$. This new vector now
fills in the remaining $n-1$ positions in $\hat{\bf X}$.

More formally, the construction is as follows. After generating ${\bf X}$ and
$(\hat{X}'_1, \hat{X}''_1)$ independently from their respective distributions, we define
$$
R=|\{ X_2,\ldots,X_n\} \cap \{\hat{X}'_1, \hat{X}''_1\}|.
$$
There are three cases.

\vspace{0.25in}

\noindent Case 0: $R = 0$.
In this case, set $(\hat{X}'_1, \hat{X}''_1, \hat{X}_2, \ldots, \hat{X}_n) =
(\hat{X}'_1, \hat{X}''_1, X_2, \ldots, X_n)$.

\vspace{0.10in}

\noindent Case 1: $R = 1$.
If say, $\hat{X}'_1$ equals $X_J$, then set $\hat{X}_i = X_i$ for $2 \le i \le n, i \neq J$
and let $\hat{X}_J$ be drawn uniformly from ${\cal A} \setminus
\{\hat{X}'_1, \hat{X}''_1, X_2, \ldots, X_n\}$.

\vspace{0.10in}

\noindent Case 2:  $R =2$.  If $\hat{X}'_1 = X_J$ and $\hat{X}''_1 = X_K$, say, then set $\hat{X}_i =
X_i$ for $2 \le i \le n, i \notin \{J, K \}$, and let $\{\hat{X}_J, \hat{X}_K\}$ be a simple
random sample of size 2 drawn from  ${\cal A} \setminus \{\hat{X}'_1,
\hat{X}''_1, X_2,\ldots, X_n\}$.

\vspace{0.10in}

The following proposition follows from Proposition \ref{zerosrs}, the
representation of the distribution (\ref{(**)}) as the product (\ref{factor}),
and that fact that conditional on $\{{\hat X}_1',{\hat X}_1'' \}$,
the above construction leads  to sampling uniformly by rejection from ${\cal A}\setminus \{{\hat X}_1',{\hat X}_1'' \}$.

\begin{proposition}
\label{wwscopr}
Let ${\bf X}=(X_1,\ldots,X_n)$ be a simple random sample of
size $n$ from ${\cal A}$ and let $(\hat{X}'_1, \hat{X}''_2) \sim
q(u,v)$ be independent of ${\bf X}$. If $\hat{X}_2,\ldots,\hat{X}_n$,
given $\hat{X}'_1, \hat{X}''_1, X_2, \ldots, X_n$, are constructed as above,
 then  $(\hat{X}'_1, \hat{X}''_1, \hat{X}_2, \ldots,
\hat{X}_n)$ has distribution (\ref{(**)}), and with $U$ an independent uniform
variate on $[0,1]$,
\beas
W^* &=&  U \hat{X}'_1 + (1-U) \hat{X}''_1 + \hat{X}_2 + \cdots + \hat{X}_n\\
W &=& X_1 + \cdots + X_n
\enas
\noindent is a realization of $(W, W^*)$ on a joint space where $W^*$ has the $W-$zero
biased distribution.
\end{proposition}

Under the moment conditions in (\ref{moment}), we have now the ingredients
to show that a bound of order $1/n$ holds, for smooth functions, for the
normal approximation of $W=\sum_{i=1}^n X_i$. First, define
\bea
\nonumber
\langle k \rangle &=& \sum_{a \in {\cal A}}a^k, \\
\label{C_1}
C_1(N,n,{\cal A})&=& \sqrt{8}\Big( \frac{\sigma^2}{4n^2}+ \langle 6 \rangle \alpha^2 + \beta^2 + \gamma^2 (n-1)^2 + \eta^2  \Big)^{1/2} \nonumber \\
\mbox{and}\\
\label{C_2} C_2(N,{\cal A})&=&11{\langle 4 \rangle}+\frac{45}{N},
\ena
where $\alpha, \beta,  \gamma$ and $\eta$ are given in (\ref{alpha}), (\ref{beta}), (\ref{gamma}), and (\ref{eta}) respectively.

\begin{theorem} \label{srs} Let $X_1, \ldots, X_n$ be a simple random sample of size $n$ from a set of $N$ real numbers ${\cal A} $ satisfying (\ref{moment}). Then
with $W=\sum_{i=1}^n X_i$, for all $h$ with four bounded derivatives we have
\begin{equation}
\label{thmbnd}
|Eh(W/\sigma)-\Phi h|  \leq \frac{1}{3  \sigma}C_1(N,n,{\cal A})||h^{(3)}||+ \frac{1}{8  \sigma^2}C_2(N,{\cal A})||h^{(4)}||.
\end{equation}
Further, if $n \rightarrow \infty$ so that $n / N \rightarrow f \in (0,1)$,
then it follows that
$$
|Eh(W/\sigma)-\Phi h|  \leq n^{-1} \{ B_1
  ||h^{(3)}|| + B_2  ||h^{(4)}|| \} (1 + o(1)),
$$
where
$$
 B_1 = \frac{\sqrt{8}}{3}\Big( \frac{f(1-f)}{4}+n^2 \langle 6 \rangle + 2 \left( \frac{f}{1-f} \right)^2  \Big)^{1/2}
       \big(f(1-f)\big)^{-1/2}
\mbox{ and } B_2 = \frac{1}{8} \big( 11 n \langle 4 \rangle + 45 f \big)\big(f(1-f)\big)^{-1} .
$$
\end{theorem}

We see as follows that this bound yields a rate $n^{-1}$ quite generally when values in $\cal A$ are ``comparable.'' For example, suppose that $Y_1,Y_2, \ldots$ are independent copies of a
nontrivial random variable $Y$ with $EY^6< \infty$ and $EY^2=1$. If $N$ is say, even, let the elements
of ${\cal A}$ be equal to the $N/2$ values
 $Y_1/(2 \sum_1^{N/2} Y_j^2)^{1/2},
\ldots, Y_{N/2}/(2 \sum_1 ^{N/2} Y_j^2)^{1/2}$ and their negatives. Then,
this collection satisfies (\ref{moment}) and (\ref{<2>}), and by the law of large numbers, a.s. as
$N \rightarrow \infty$, the terms $n \langle 4 \rangle$ and $n^2 \langle 6 \rangle$ converge to constants. Specifically,
$$
n \langle 4 \rangle \rightarrow f E Y^4  \quad \mbox{and} \quad
n^2 \langle 6 \rangle \rightarrow f^2 E Y^6,
$$
and so $B_1$ and $B_2$ are asymptotically constant, and the
rate $1/n$ is achieved over the class of functions $h$ with
bounded derivatives up to fourth order.

\noindent {\bf Proof.}  Both $Eh(W)$ and the upper bound in (\ref{thmbnd})
are continuous functions of $\{a_1,\ldots,a_N\}$. Hence, since any collection
of $N$ numbers ${\cal A}$ is arbitrarily close to a collection of $N$
distinct numbers, it suffices to prove the theorem under the assumption
that the elements of ${\cal A}$ are distinct.

We apply Theorem \ref{zerogen}.
Constructing $W$ and $W^*$ as in Proposition \ref{wwscopr}, and using
standard inequalities and routine computations, one can show that
\beas
\mbox{Var}\{E(W^*-W|W)\} \leq 2
\{\mbox{Var}\Big(\frac{1}{n}W\Big)+\mbox{Var}(A)\}
\enas
where
\beas
A&=&\alpha \sum\limits_{x \in {\bf Y}} x^3 +  \beta \sum\limits_{x \in {\bf Y}} x
\sum\limits_{x \in {\bf Y}} x^2
+ \gamma (\sum\limits_{x\in {\bf Y}} x)^3 + \eta \sum\limits_{x \in {\bf Y}} x,
\enas
and
\bea \label{alpha}
\alpha&=&\frac{n-1}{N(N-n)}-1
\\
\label{beta}
\beta&=&\frac{-2(n-1)}{N(N-n+1)}+\frac{n-3}{N(N-n)}-\frac{1}{N}
\\
\label{gamma}
\gamma &=& \frac{-2}{N (N-n)(N-n+1)}
\\
\label{eta}
\eta&=&\frac{-N+3}{N(N-n)},
\ena
and that
\beas
\mbox{Var}(E\{W^*-W| W \} ) \le C_1^2(N,n,{\cal A}),
\enas
where $C_1(N,n,{\cal A})$ is given in (\ref{C_1}). Again, straightforward computations give that
\beas
E(W^*-W)^2 &\leq& 11 \sum_{a \in {\cal A}} a^4 + \frac{45}{N}=C_2(N,{\cal A}).
\enas
Details can be found in the technical report (\ref{grtech}).

\section*{Bibliography}
\begin{enumerate}

\item \label{brs}  Baldi, P. Rinott, Y.  and Stein, C.  (1989). {A normal approximations for  the number of local maxima of a random function on a graph}, In {\sl Probability, Statistics and Mathematics, Papers in Honor of Samuel Karlin.}  T. W. Anderson, K.B. Athreya and D. L. Iglehart eds., Academic Press, 59--81.

\item \label{ADB90}  Barbour, A.D. (1990). {Stein's method for diffusion approximations.}
{\sl Probab. Theory Related Fields } {\bf 84}, 297--322.


\item \label{BG93} Bolthausen, E. and G\"otze, F. (1993). The rate of convergence for multivariate sampling statistics. {\sl Ann. Statist.} {\bf 21}, 1692--1710.


\item \label{cac2} Cacoullos, T., Papathanasiou, V. and Utev, S. (1994). Variational inequalities with examples and an application to the central limit
theorem.  {\sl Ann. Probab.} {\bf 22}, 1607--1618.


\item \label{yosefication} Dembo, A. and Rinott, Y. (1996). Some examples of Normal approximations by Stein's method. In {\sl Random Discrete Structures}, IMA volume 76, 25--44.
Aldous, D. and Pemantle, R. Eds., Springer-Verlag.

\item \label{Feller2} Feller, W. (1971). An Introduction to Probability and Its Applications. Volume 2. 2$^{nd}$ edition. John Wiley and Sons.

\item \label{foll} Folland, G. (1984) Real Analysis: Modern
Techniques and Their Applications. John Wiley and Sons.

\item \label{grtech} Goldstein, L. and Reinert, G. (1997) {Stein's Method and the Zero Bias Transformation with Application to Simple Random Sampling:
Technical Report}

\item \label{gr96} Goldstein, L. and Rinott, Y. (1996). {On
multivariate normal approximations by Stein's method and size bias
couplings. {\sl J. Appl. Prob.} {\bf 33}, 1--17}.
arXiv:math.PR/0510586

\item \label{g91} G\"otze, F. (1991). {On the rate of convergence in the multivariate CLT}. {\sl Ann. Probab.}  {\bf 19}, 724--739.

\item \label{hoch78} Ho, S.-T. and Chen. L.H.Y. (1978). An $L_p$ bound for the remainder in a combinatorial central limit theorem. {\sl Ann. Statist.} {\bf 6}, 231--249.


\item \label{yos94} Rinott, Y. (1994). On normal approximation rates
for certain sums of dependent random variables.  {\sl J. Computational and Applied Math.}
{\bf 55}, 135--143.

\item \label{rinrot}  Rinott, Y. and Rotar, V. (1996): {A multivariate CLT for local dependence with $n^{-1/2} \log n$ rate and
applications to multivariate graph related statistics.} {\sl J. Multivariate Analysis}. {\bf 56}, 333-350.

\item \label{st72} Stein, C. (1972). {A bound for the error in the normal approximation to
the distribution of a sum of dependent random variables.} {\sl Proc. Sixth Berkeley Symp. Math. Statist. Probab.} {\bf 2}, 583--602. Univ. California Press, Berkeley.

\item \label{st86} Stein, C. (1986).  {\sl Approximate Computation of Expectations}.   IMS, Hayward, California.

\item \label{st92} Stein, C. (1992). {A way of using auxiliary randomization}. In {\sl  Probability Theory,}  159--180. Walter de Gruyter, Berlin - New York.

\end{enumerate}
\end{document}